\documentclass[11pt]{article}
\newcommand{\comm}[1]{}

\usepackage{amssymb}\usepackage{graphicx}
\usepackage{amsmath}
\usepackage{times}\hfuzz=10pt 
\usepackage[numbers]{natbib}
\usepackage{epsfig}
\def\citet{\cite}

\newtheorem{theorem}{Theorem}

\def\e{\varepsilon}

\def\defi{\stackrel{{\scriptscriptstyle \Delta}}{=}}

\def\a{\alpha}
\def\d{\delta}
\def\o{\omega}
\def\O{\Omega}

\def\F{{\cal F}}
\def\w{\widehat}
\def\Ind{{\,\rm Ind\,}}
\def\Ind{{\mathbb{I}}}

\def\Re{{\rm Re\,}}
\def\Im{{\rm Im\,}}

\def\R{{\bf R}}
\def\E{{\bf E}}

\def\Z{{\cal Z}}

\def\b{\beta}
\def\s{\delta}
\def\g{\gamma}
\def\C{{\bf C}}

\def\ww{\widetilde}

\def\t{\theta}
\def\oo{\bar}
\def\s{\sigma}

\newcommand{\be}{\begin{equation}}
\newcommand{\ee}{\end{equation}}
\newcommand{\bd}{\begin{displaymath}}
\newcommand{\ed}{\end{displaymath}}
\newcommand{\ba}{\begin{array}{ll}}
\newcommand{\ea}{\end{array}}
\newcommand{\baa}{\begin{eqnarray}}
\newcommand{\eaa}{\end{eqnarray}}
\newcommand{\baaa}{\begin{eqnarray*}}
\newcommand{\eaaa}{\end{eqnarray*}}
\font\sm=cmr10

\def\KH{{\scriptscriptstyle K\!H}}
\def\KK{{\scriptscriptstyle K}}

\def\T{{\cal T}}

\def\oo{\bar}
\def\aa{a}

\def\a{\alpha}

\def\hh{h}
\def\HH{H}
\def\a{a}
\def\b{b}
\def\ew{\left(e^{i\o}\right)}
\def\T{{\mathbb{T}}}
\def\ZZ{{\mathbb{Z}}}

\def\HHH{{\rm H}}
\def\ZZ{{\mathbb{Z}}}
\def\muu{\theta}
\def\subitem{\item}
\date{Submitted November 27, 2014; revised July 4, 2015}
\title{Sub-ideal causal smoothing filters for the real sequences}
\author{
Nikolai Dokuchaev\\  {\sm Department of Mathematics \& Statistics,
Curtin University,}\\
{\sm  GPO Box U1987, Perth, 6845 Western Australia}}
 
\begin{document}
\def\break{}%
\def\brea{}
\def\breakk{}
\def\breakSP{\nonumber\\&&}
\def\breaSP{\nonumber\\}
      \maketitle
\begin{abstract}  The paper considers  causal smoothing of the
real sequences, i.e.,
discrete time processes in a deterministic setting. A family of  causal linear time-invariant
filters is suggested.  These filters approximate the gain decay for some non-causal
ideal smoothing
filters with transfer functions vanishing at a point of the unit circle and such
that they transfer processes into predictable ones.
In this sense, the suggested filters are  near-ideal;
a faster gain decay   would lead to the loss of causality.   Applications to predicting algorithms are discussed
and illustrated by experiments with  forecasting of  autoregressions  with the coefficients that are deemed to be untraceable.
\end{abstract}
\index{HIGHLIGHTS (FOR ONLINE VERSION)A family of smoothing causal
linear time-invariant filters for continuous time processes is
suggested such that the frequency response approximates the real
unity uniformly on an arbitrarily large interval. > In addition,
almost exponential damping of the energy on the higher frequencies
is ensured. > These filters are near-ideal meaning that a faster
decay of the frequency response would lead to the loss of causality}
\maketitle
\section{Introduction}  The
paper studies causal smoothing of the  discrete time processes.
For many applications, it is preferable to replace a process by  a more smooth process.
In continuous time setting, smoothness is associated with predictability. Smooth  analytic functions are predictable, i.e., their values on any interval define uniquely their values outside of this interval, and an ideal low-pass filter converts a function into an analytic one.
For discrete time processes,  it is not obvious how to
define an analog of the continuous time analyticity and smoothness. A classical approach is to consider  predictability instead of analyticity.  So far, the predictability criterion  for stochastic
Gaussian stationary discrete time processes in the frequency domain setting are
given by the classical Szeg\"o-Kolmogorov Theorem. This theorem says that
the optimal prediction error is zero  if  \baa \int_{-\pi}^\pi \log\phi\ew d\o=-\infty, \eaa where $\phi$ is the spectral density; see Kolmogorov \citet{Ko},
Szeg\"o \citet{Sz,Sz1}, Verblunsky \citet{Ver}, and more recent literature reviews in
\citet{Bin,Si}. This means that a stationary Gaussian  process  is
predictable if its spectral density is vanishing on a part of the
unit circle $\{z\in\C:\ |z|=1\}$, i.e., if the process is
"band-limited" in this sense. This result was expanded on more
general stable stochastic processes   allowing spectral representations with
spectral density via processes with independent increments; see, e.g., \citet{CS}.

The stochastic setting is the most common in
causal smoothing and sampling; see, e.g.,
\citet{AU,rema,B,Br,Bin,D08,D10e,D12a,D12sp,F94,PFG,PFG1,FKR,Higgins,Ko,Mi,Mizes}.
A transition to the predicability of discrete time processes in deterministic setting   is non-trivial
 and related to the concept of the randomness  for the real sequences in the pathwise setting without a probability measure.
 There are many classical works devoted to this  important concept, starting from Mises  \cite{Mizes}, Church  \cite{Chu}, Kolmogorov
 \cite{Ko}, Loveland \cite{Lo}; see the references in \citet{Li}.
 \par
 It was found that real sequences are predictable  if their  Z-transform vanishes on an arc  of the unit circle \citet{D12sp} on the complex plane or at a point $z=-1$ of the unit circle \citet{D12a}.
Therefore, smoothing can be interpreted as reduction of the energy on
the higher frequencies. In particular, an ideal low-pass filter is a
smoothing filter. This filter is non-causal, i.e., it
requires the future value of the process. Similarly,  a filter with too high rate of decay  of the frequency response at a certain point of the unit circle
 also cannot be causal, since causality is inconsistent with predictability of  outputs described in \citet{D12a}.
\par
The present paper readdresses the problem of causal smoothing of the discrete time processes in the deterministic pathwise setting, without probabilistic assumptions.
We suggest a family of causal smoothing filters
that can be arbitrarily close to some ideal non-causal smoothing filters  defined by equation (\ref{NC}) below.
 The suggested filters are near-ideal in the sense that they ensure "almost"
ideal  rate of damping the energy at the point $z=-1$; a faster decay of the
frequency response is impossible for causal filters.  This follows from
predictability criterion \citet{D12sp}. In fact, the particular
reference family    of  non-causal ideal filters (\ref{NC}) was selected because these
filters
 transfer  non-predictable processes into predictable ones satisfying the criterion from  \citet{D12a}.
 Similar approach was used in \cite{D10f} for the continuous time setting.
\par
The suggested near-ideal filters are  discrete time causal linear
time-invariant filters (LTI filters); they are
represented as convolution integrals over the historical data and  approximate the
real unity uniformly on an arbitrarily large part of the unit circle.
\par
 It appears that these causal filters can be used to improve the performance of  predictors  suggested in \citet{D12a}.
  These are robust with respect to the noise contamination, and the error caused by
   a high-frequency noise depends on the intensity of this noise.  A near-ideal smoothing  filter cannot remove   the high-frequency noise entirely but still can reduce it. This approach is discussed in  Section \ref{secF}, where we present some numerical experiments  with the suggested near-ideal filters
applied to forecasting of  autoregressions in a setting where the autoregression coefficients are deemed to be untraceable.
\subsubsection*{Some definitions and notations}
We denote by $\ZZ$  the set of all integers.
\par
For $r\in[1,+\infty]$, we denote by $\ell_r$ the set of all
sequences $x=\{x(t)\}_{t\in\ZZ}\subset\R$, such that
$\|x\|_{\ell_r}=\left(\sum_{t=-\infty}^{\infty}|x(t)|^r\right)^{1/r}<+\infty$
for $r\in[1,\infty)$ or  $\|x\|_{\ell_\infty}=\sup_t|x(t)|<+\infty$
for $r=+\infty$. We denote by $\ell_r^+$ the set of all sequences
$x\in\ell_r$ such that $x(t)=0$ for $t<0$.
\par
We denote by $L_r(-\pi,\pi)$ the usual Banach space of complex valued $L_r$-integrable functions $x:[-\pi,\pi]\to\C$.

Let $D^c\defi\{z\in\C: |z|> 1\}$, and let $\T=\{z\in\C: |z|=1\}$.
\par
For  $x\in \ell_1$ or $x\in \ell_2$, we denote by $X=\Z x$ the
Z-transform  \baaa X(z)=\sum_{t=-\infty}^{\infty}x(t)z^{-t},\quad
z\in\C. \eaaa Respectively, the inverse Z-transform  $x=\Z^{-1}X$ is
defined as \baaa x(t)=\frac{1}{2\pi}\int_{-\pi}^\pi
X\left(e^{i\o}\right) e^{i\o t}d\o, \quad t=0,\pm 1,\pm 2,....\eaaa
 If $x\in \ell_2$, then $X|_\T$ is defined as an element of
$L_2(\T)$, i.e., $X\ew\in L_2(-\pi,\pi)$. If $X\ew\in L_1(-\pi,\pi)$, then $x=\Z^{-1}X$
is defined as an element of $\ell_\infty$.

 Let  $\HHH^2(D^c)$ be the Hardy space of functions that are holomorphic on
$D^c$ including the point at infinity  with finite norm
$\|h\|_{\HHH^2(D^c)}=\sup_{\rho>1}\|h(\rho e^{i\o}))\|_{L_2(-\pi,\pi)}$.
Note that Z-transform defines a bijection
between the sequences from $\ell_2^+$ and the restrictions (i.e.,
traces) $X|_{\T}$ of the functions from $\HH^2(D^c)$ such that  $\overline{X\ew} =X\left(e^{-i\o}\right)$ for $\o\in\R$; see, e.g., \cite{Linq}, Section 4.3.
If $X\ew\in L_1(-\pi,\pi)$ and $\overline{X\ew} =X\left(e^{-i\o}\right)$, then $x=\Z^{-1}X$
is defined as an element of $\ell_\infty^+$.
\section{Problem setting}
Let $x(t)$ be a discrete time process, $t\in\ZZ$. The output of a
linear filter is the process \baaa y(t)=\sum_{s=-\infty}^{\infty}
\hh(t-s)x(s),\eaaa where $\hh:\ZZ\to\R$ is a given impulse
response function.

If $\hh(t)=0$ for $t<0$, then the output of the corresponding filter
is \baaa y(t)=\sum_{s=-\infty}^{t} \hh(t-s)x(s). \eaaa In
this case, the filter and the impulse response function are said to
be causal. The output of a causal  filter at time $t$ can be
calculated using only past historical values $x(s)|_{s\le t}$
of the currently observable  input process.

The goal is to approximate $x$ by a "smooth" filtered process
$y$ via selection of an appropriate causal impulse response
function $\hh$.
\par
We are looking for families of the causal smoothing impulse response
functions $\hh$  satisfying the following conditions.
\begin{itemize}
\item[(A)] The  outputs $y$ approximate inputs $x$; an arbitrarily close approximation
can be achieved via selection of a  filter from this family.
\item[(B)] The spectrum of the output $y$ vanishes on higher frequencies.
\index{\item[(C)] The effectiveness of this family in the damping of the
higher frequencies cannot be exceeded; any faster decay of the
frequency response  would lead to the loss of causality.}
\item[(C)] The effectiveness  of the damping on the energy on the
higher frequencies  approximates the effectiveness of some reference
family of non-causal smoothing  filters  that transfer processes into predictable ones, i.e., such that the future values are uniquely
defined by the past values.
\end{itemize}
\par
Note that it is not a trivial task to satisfy Conditions (A)-(C) simultaneously. For example,
 there are sets of
 ideal low-pass filters such that the distance of these sets  from the set of
all causal filters is zero.  In \cite{rema}, this was shown  for the set of low-pass filters with increasing
pass interval $[-\Delta,\Delta]$, where $\Delta\in(0,\pi)$. However, Condition (B) is not satisfied
for the corresponding causal approximations. Moreover,
 all known causal approximations of ideal filters do not feature
 zero values for the transfer functions. The present paper suggests transfer functions vanishing at $z=-1$
together with their derivatives.
\subsubsection*{The targeted properties of the near ideal filters}
Our purpose is to   construct  a family of causal
filters  such  that the Conditions (A)--(C) are satisfied.
We will be using   a reference family of "ideal" smoothing filters with the
frequency response \baa M_{\muu ,q}\ew=\exp\left(-\frac{\muu }{|1+e^{i\o}|^q} \right),
\label{NC}\eaa where  $q>1$ and $\muu>0$ are parameters. For these filters, Condition (A) is satisfied  as
$\muu \to 0$, and Conditions (B) is satisfied for all $\muu >0$.
However, these filters are non-causal because, for any $x\in
\ell_2$, the values $x(t+1)$ of the output processes of these filters are weakly
predictable at time $t$ \cite{D12a}. This is since  $M_{\muu ,q}\ew \to 0$
fast enough
as $\o\to -\pi$.
\par
For a given integer $m\ge 1$, we will construct a family of causal filters with impulse responses
$\hh_{\aa}\in \ell_\infty^+$ and with the
corresponding Z-transforms $\HH_{\aa}=\Z h_\aa$, where $\aa\in(0,1)$ is a parameter.
For this family, the
following more special Conditions (a)-(c) will be satisfied (the number $m$ is used in Condition (b1) below).
\begin{itemize}
\item[(a)]  {\em  Approximation of the identity operator:}
\subitem(a1)  $\sup_{\o\in[0,\pi], \aa}|\HH_{\aa}\ew|<+\infty$.
\subitem(a2) For any $\O>0$, $\HH_{\aa}\ew\to 1$ as $\aa\to 1-0$
uniformly in $\o\in[-\O,\O]$.
\subitem(a3) For any  $X\ew\in L_1(-\pi,\pi)$ and any $x=\Z^{-1}X$, \baaa \|y_{\aa}(\cdot)-x(\cdot)\|_{\ell_\infty}\to 0\quad
\hbox{as}\quad \aa\to 1-0,\eaaa where     $y_{\aa}$ is the output
process \baaa y_{\aa}(t)=\sum_{\tau=-\infty}^{t}
\hh_{\aa}(t-\tau)x(\tau).\eaaa
\item[(b)] {\em The spectrum is vanishing at a point at $\T$:} \subitem(b1) For all $\aa$, $\HH_{\aa}\ew$ is $m$ times differentiable at $\o=\pi$, and
\baaa
\HH_{\aa}(-1)=0,\qquad \frac{d^k H_\aa}{d\o^k}\ew\Bigl|_{\o=\pi}=0,\quad k=1,...,m.
\eaaa
\subitem(b2) For any $\e>0$, there exists $\d>0$ and $\aa\in (0,1)$ such that
 \baaa
\sup_{\o\in[\pi-\d,\pi+\d]}|\HH_{\aa}\ew|<\e.\eaaa

\index{\item[(c)]  {\em Sub-ideal smoothing:} There exists $\oo\aa>0$ such that
for any $\aa>\oo\aa$  \baa \int_{-\pi}^{\pi}|\log|\HH_{\aa}\ew||d\o\to -\infty\quad \hbox{as}\quad \aa\to +\infty.
\label{delta} \eaa}
\item[(c)]  {\em Approximation  of non-causal filters
(\ref{NC}) with respect to the effectiveness in damping:} For any
$\e>0$ and any $\O\in(0,\pi)$,  $\O_0\in(\O,\pi)$,  $\O_1\in(\O_0,\pi)$
there exists  $\muu >0$, $q>1$, $\aa>0$ such that  \baa &&|\HH_{\aa}\ew-1|\le \e,\quad \o\in[-\O,\O],\label{un1}\\
&&|\HH_{\aa}\ew|\le |M_{\muu ,q}\ew|,\breakk \quad\quad\o\in[-\O_1,-\O_0]\cup [\O_0,\O_1].
\label{un2} \eaa
\index{\baaa
|\HH_{\aa}\ew|\le |M_{\muu ,q}\ew|+ \e,\quad \o\in[-\pi,\pi]\backslash[-\O_0,\O_0]. \eaaa}
\end{itemize}
\par
Conditions (a)-(c) represent particular versions of less specific Conditions (A)-(C).
In particular, estimate (\ref{un2}) ensures that  Condition (C) is satisfied.
\section{A family of near-ideal smoothing filters}
 Let a real number $p\in (1/2,1)$ and integers $N\ge 1$ and $m\ge 1$ be given.
 For the real numbers $\a\in(0,1)$,  we define transfer functions  \baa
\HH_{\aa}(z)=\left(\exp\frac{(1-a)^p}{z+\a}+G_\aa(z)\right)^m, \quad z\in\C,\label{K}\eaa
where
\baaa
G_{\aa}(z)=-\xi(\a,p)+\frac{\gamma(\a,p)}{N} \left((-1)^Nz^{-N}-1\right),
\label{G}\eaaa
and where \baaa
&&\xi(\a,p)=\exp[-(1-\a)^{p-1}],\def\break{}%
\quad\breakk \g(\a,p)= |1-\a|^{p-2}\xi(\a,p).  
\eaaa

\par
We consider the set $\{\HH_\aa\}_{a\in(0,1)}$ of transfer functions (\ref{K}) with  a fixed triplet $(m,N,p)$.

\begin{theorem}\label{ThM} Conditions (a)-(c) are
satisfied  for  the family of filters defined by the transfer
functions $\{\HH_{\aa}\}_{a\in(0,1)}$. (Therefore, Conditions
(A)-(C) are satisfied for this family).
\end{theorem}
\par
{\em Proof of Theorem \ref{ThM}.} Let us assume first that $m=1$.

Clearly, the functions $\HH_{\aa}$ are holomorphic in $D^c$ and bounded in $D^c\cup\T$ for any $\aa\in(0,1)$.
Hence the inverse Z-transforms
$\hh_{\aa}=\Z^{-1}\HH_{\aa}$ are causal impulse responses, i.e.,
$\hh_{\aa}(t)=0$ for $t<0$; see, e.g., \cite{Linq}, Theorem 4.3.2.
\par
Let $f(\a)=(1-\a)^p$ and $\Psi_{\aa}(z)=f(\a)(z+\a)^{-1}$.
By the definitions,  $H_{\aa}(z)=\exp\Psi_{\aa}(z)+G_{\aa}(z)$, and
\baaa
\Psi_{\aa}\ew=f(\a)\frac{\cos(\o)+\a-i\sin(\o)}{(\cos(\o)+\a)^2+\sin(\o)^2}.
\eaaa
\par
Let us prove that Condition (a) holds.

Clearly, $|G_{\aa}\ew|\to 0$ as $\aa\to 1-0$ uniformly in $\o\in(-\pi,\pi]$.

Let $\o_a\in (\pi/2,\pi)$ be such that $\cos(\o_a)+\a=0$. We
have that $\Re\Psi_{\aa}\ew>0$ for all $\o\in[-\o_a,\o_a]$ and
$\Re\Psi_{\aa}\ew<0$ for all $\o\in[-\pi,\o_a)\cup(-\o_a,\pi]$.

Further, we have that \baaa \inf_{\o\in[-\o_a,\o_a] }|e^{i\o}+a|\ge \sqrt{1-a^2}.\eaaa  Hence
\baaa \sup_{\o\in[-\o_a,\o_a] }|\Psi_{\aa}\left(e^{i\o}\right)|\le \frac{f(a)}{\sqrt{1-a^2}}=\frac{(1-a)^{p-1/2}}{(1+a)^{1/2}}\le 1.\eaaa
Therefore, the value
$|H_{\aa}\ew|$ is uniformly bounded in $\a,\o$.
Hence  Condition (a1) holds.

Further, we have that \baaa \o_a\to \pi-0\quad\hbox{as}\quad
\a\to 1. \eaaa Hence, for any $\O\in[0,\pi)$, we have that \baaa \sup_{\o\in[-\O,\O] }|\Psi_{\aa}\left(e^{i\o}\right)|\to 0\quad\hbox{as}\quad
\a\to 1\eaaa and
\baaa \sup_{\o\in[-\O,\O] }|H_{\aa}\left(e^{i\o}\right)-1|\to 0\quad\hbox{as}\quad
\a\to 1.\eaaa
Hence Condition (a2) holds.

 Let as show that Condition (a3) holds.
 Let $Y_{\aa}=\HH_{\aa} X$.  By Condition
(a2), $ Y_{\aa}\ew\to X\ew$ as $\aa\to 1-0$ for all $\o\in\R$.
Clearly, there exists  $\aa_0\in(0,1)$ and $c_0>0$ such that  $\sup_{\o,\aa\ge \aa_0}|\HH_{\aa}\ew|\le c_0$. Hence $|Y_{\aa}\ew-X\ew|\le
(c_0+1) |X\ew|$. By the assumptions, $X\ew=\Z x\in L_1(-\pi,\pi)$. By the Lebesgue
Dominance Theorem, it follows that \baaa \left\|Y_{\aa}\ew-X\ew
\right\|_{L_1(-\pi,\pi)}\to 0\quad\hbox{as}\quad \aa\to 1-0.
 \eaaa
Therefore, Condition (a3) holds and Condition (a) holds.
\par
Let us show that Condition (b) holds. We have that
\baaa
&&\exp\left(\Psi_{\aa}\left(e^{i\pi}\right)\right)
=\exp\left(f(\a)\frac{-1+\a}{(1-\a)^2}\right)\breakk
=\exp\left(-\frac{(1-\a)^p}{1-\a}\right)\breakSP=\xi(\a,p)=-G_{\aa}(-1).\eaaa
Hence $\HH_{\aa}(-1)=\exp\left(\Psi_{\aa}\left(-1\right)\right)+G_{\aa}(-1)=0$.

Let us show that $\frac{dH_{\aa}}{d\o}\ew\Bigl|_{\o=\pi}=0$.
Let  \baaa
&&r(\o)=\Re\exp\left(\Psi_{\aa}\left(e^{i\o}\right)\right),\breakk\quad s(\o)=\Im\exp\left(\Psi_{\aa}\left(e^{i\o}\right)\right),\quad
\breakSP q(\o)=\frac{1}{N}\Im\left(e^{-iN(\o-\pi)}-1\right).\eaaa
Clearly, the function $\Psi_{\aa}\ew$ is differentiable in $\o\in \R$ for any $\aa$, as well as functions $r(\o)$, $s(\o)$, and $q(\o)$.
 In addition, we have  that
$r(\o)=r(\pi-\o)$. Hence $r(\o)$
is even about the point $\o=\pi$ and  differentiable. This implies that
\baaa
\frac{dr}{d\o}(\o)\Bigl|_{\o=\pi}=0.
\eaaa
 By the definitions,  $s(\o)=\exp(\Re\Psi_{\aa})\sin(\Im\Re\Psi_{\aa})$ and
$\exp(\Re\Psi_{\aa})\ew  \to \xi(\a,p)$ as $\o\to\pi$. We have that $s(\pi)=q(\pi)=0$.  The L'H\^opital's rule gives  that
\baaa
&&\lim_{\o\to\pi}\frac{\frac{ds(\o)}{d\o}}{\frac{dq(\o)}{d\o}}= \lim_{\o\to\pi}\frac{s(\o)}{q(\o)} \breakk =
\lim_{\o\to\pi}\frac{\xi(\a,p)\sin\left(\frac{-(1-\a)^p\sin(\o)}{(a+\cos(\o))^2+\sin(\o)^2}\right)}{\frac{1}{N}\sin(N(\o-\pi))}  \\&& =-\xi(\a,p)\frac{(1-\a)^p}{(a-1)^2}=-\g(\a,p).
\eaaa
Clearly, $\frac{dq(\o)}{d\o}|_{\o=\pi}=-1$.  Hence
\baaa
\frac{d}{d\o}\exp\Psi_{\aa}\left(e^{i\o}\right)\Bigl|_{\o=\pi}=i\frac{ds(\o)}{d\o}|_{\o=\pi}=i\g(\a,p).
\eaaa
On the other hand, \baaa
\frac{dG_{\aa}\ew}{d\o}\Bigl|_{\o=\pi} \breakk =\frac{d}{d\o}\left(-\xi(a,p)+\frac{\g(\a,p)}{N}\left(e^{-iN(\o-\pi)}-1\right)\right)\Bigl|_{\o=\pi}\breaSP=-i\g(a,p).
\eaaa
Hence
\baaa
&&\frac{d}{d\o}\HH_{\aa}\left(e^{i\pi}\right)\Bigl|_{\o=\pi}\breakk=\frac{d}{d\o}\exp\left(\Psi_{\aa}\left(e^{i\pi}\right)\right)\Bigl|_{\o=\pi}
+\frac{dG_{\aa}\ew}{d\o}\Bigl|_{\o=\pi} =0.
\eaaa
 Therefore, Condition (b1) holds.
\par
Let us show that Condition (b2) holds. We have that \baaa
\Re\Psi_{\aa}\ew=\frac{f(\a)(\cos(\o)+\a)}{(\cos(\o)+\a)^2+\sin(\o)^2}=f(a)\frac{\Re(e^{i\o}+a)}{|e^{i\o}+a|^2}.
\eaaa
We have that $-\Re(e^{i\o}+a)/|e^{i\o}+a|$ is non-decreasing in $\o\in[\o_a,\pi]$ and converges to $1$ as  $\o\to\pi-0$, and $1/|e^{i\o}+a|$  is non-decreasing in $\o\in[\o_a,\pi]$ and converges to $(1-a)^{-1}$ as  $\o\to\pi-0$. Hence the product of these functions,
$-f(a)^{-1}\Re\Psi_{\aa}\ew$, is non-decreasing in $\o\in[\o_a,\pi]$ . Hence we can select $\w\o_{\aa}\in[\o_a,\pi]$ such that $-\Re\Psi_{\aa}\ew\ge  -\Re\Psi_{\aa}\left(e^{i\pi}\right)/2$ for all $\o\in[\w\o_{\aa},\pi]$, i.e.,  $\Re\Psi_{\aa}\ew\le  \Re\Psi_{\aa}\left(e^{i\pi}\right)/2$ for all $\o\in[\w\o_{\aa},\pi]$. In addition, $\Re\Psi_{\aa}\ew=\Re\Psi_{\aa}\left(e^{-i\o}\right)$. Hence
 $\Re\Psi_{\aa}\ew\le \Re\Psi_{\aa}\left(e^{i\pi}\right)/2$ for all $\o\in[\w\o_{\aa},2\pi-\w\o_{\aa}]=[\pi-\d_\aa,\pi+\d_\aa]$, where  $\d_\aa=\pi-\w\o_\aa$.

 Further,  we have that
\baaa\Re\Psi_{\aa}\left(e^{i\pi}\right)= (1-a)^p\frac{-1+a}{(-1+a)^2}\to-\infty
\quad\hbox{as}\quad \a\to 1.\eaaa
\par
 For  a given $\e>0$, let us select $\oo\aa$ such that $\Re\Psi_{\aa}\left(e^{i\pi}\right)/2<\log(\e/2)$ for all $\aa\ge\oo\aa$. In addition, we can select $\ww\aa\ge \oo\aa$ such that $|G_{\aa}\ew|\le \e/2$ for all $\aa\ge\ww\aa$ and all $\o$.
Then Condition (b2) holds with  $a=\ww\aa$ and $\d=\d_{\ww\aa}$ selected for given $\e$. Therefore, Condition (b) holds.

Let us show that Condition (c) holds.  It follows from the proof of (a) above, that, for a given $\e>0$ and $\O$, we can select $\oo\aa$ such that (\ref{un1}) holds for $\aa\ge\oo\aa$. Further, let  $q>1$ be any. For any $\O_0>\O$ and $\O_1>\O_0$,
\baaa \sup_{ \o\in[-\O_1,-\O_0]\cup [\O_0,\O_1]} |M_{\muu ,q}\ew-1|\to 0\quad \hbox{as}\quad \muu \to 0.
\eaaa
Clearly, (\ref{un2}) holds for small enough $\muu $.
Hence Condition (c) holds.

We have proved the theorem for the case where $m=1$. The extension on the case where $m>1$ is straightforward.   This completes
the proof of Theorem \ref{ThM}. $\Box$

\subsection*{Illustrative examples}
\begin{figure}[ht]
\centerline{\psfig{figure=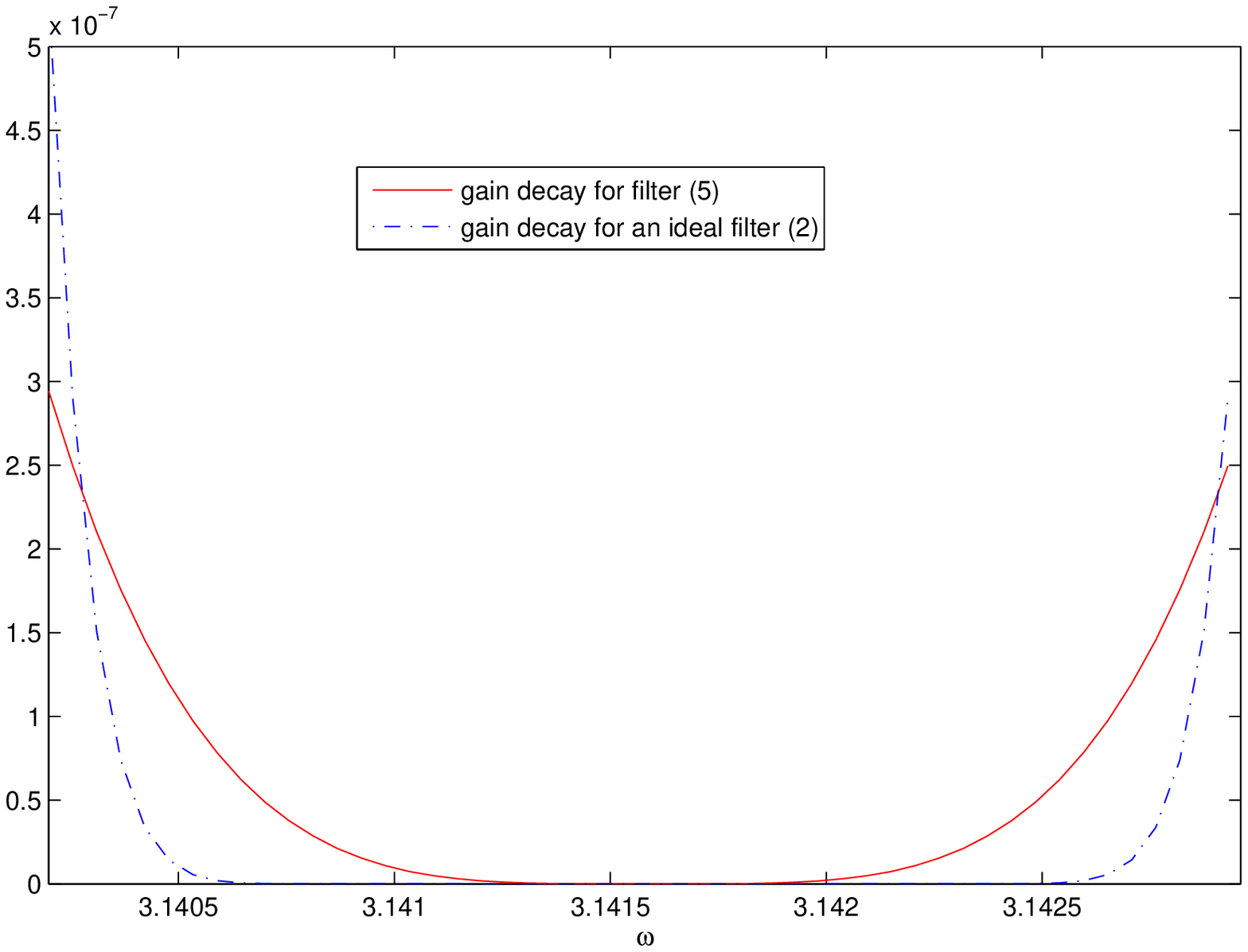,width=9cm,height=6.5cm}} \caption[]{\sm
Gain decay: the values  $|M_{\muu,q }\ew|$  for a non-causal
filter  (\ref{NC})  with $\muu =0.02$ and $q=1.01$, and $|\HH_{\aa}\ew|$ for a causal filter
(\ref{K}) with  $\a=0.99$, $p=0.6$, $N=50$, $m=2$.}
\vspace{0cm}\label{fig0}\end{figure}

\begin{figure}[ht]
\centerline{\psfig{figure=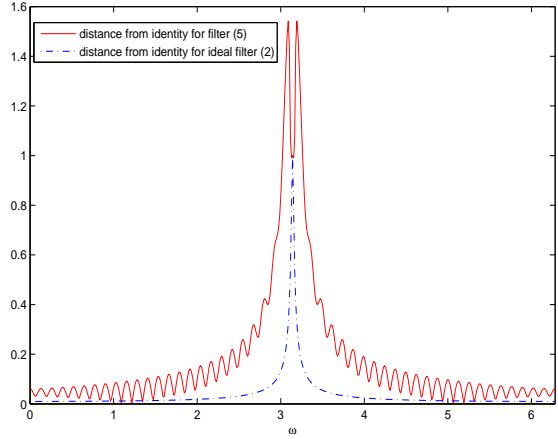,width=9cm,height=6.5cm}}
\vspace{0.5cm} \caption[]{\sm Approximation of identity operator:
shapes of the distances from 1, i.e, for the values  $|M_{\muu ,q}\ew-1|$  and $|\HH_{\aa}\ew-1|$, for a non-causal
filter  (\ref{NC})  with $\muu =0.02$ and $q=1.01$, and $|\HH_{\aa}\ew|$ for a causal filter
(\ref{K}) with  $\a=0.99$, $p=0.6$, $N=50$, $m=2$.}
\vspace{0cm}
\label{fig1}
\end{figure}

\begin{figure}[ht]
\centerline{\psfig{figure=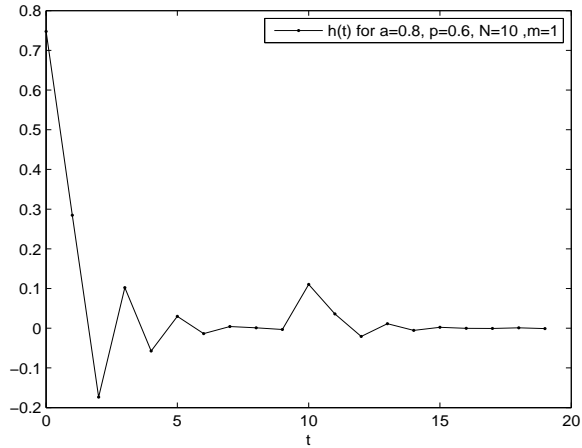, width=9cm, height=6.5cm}}
\vspace{0.5cm} \caption[]{\sm Impulse response
$h_{\a}(t)=(\F^{-1}H_{\a})(t)$ for causal filter (\ref{K})
with  $\a=0.8$, $p=0.6$, $N=10$, $m=1$.} \vspace{0cm}
\label{figh}
\end{figure}

  Figures
\ref{fig0}-\ref{figh} shows examples  of the frequency responses and the impulse  functions for the
filters described above.
\par
 Figure
\ref{fig0} shows the shapes of gain curves
$|M_{\muu ,q}\ew$  for reference non-causal filter (\ref{NC})
with $\muu =0.02$, $q=1.01$, and  $|\HH_{\aa}\ew |$ for near-ideal causal
filters (\ref{K}) with $\a=0.99$, $p=0.6$,  $N=50$, $m=2$.
\par
 Figure \ref{fig1} shows the shapes of error curves for approximation of
 identity operator on low frequencies. More precisely, it shows
 $|M_{\muu ,q}\ew -1|$ for reference non-causal filter  (\ref{NC})
and $|\HH_{\aa}\ew -1|$  for near-ideal causal
filters (\ref{K}), with the same parameters as for Figure \ref{fig0}.
\par
Figure \ref{figh} shows an example of impulse response
$h=\Z^{-1}\HH_{\aa}$ calculated  as the inverse
Z-transform for causal filter (\ref{K}) with $a=0.8$, $p=0.6$, $N=10$, $m=1$.  Since the properties of $H_{\aa}$ guarantee that $\Im h_{\aa}(t)=0$ for all $t$ and that $ h_{\aa}(t)=0$ for all $t<0$, we show the values
for  $t\ge 0$ only.
\par
In particular, these examples show  that the impulse response functions $h_\aa$ can take negative values; i.e.,  these filters do not represent
an averaging with a positive kernels.

\section{Applications to the forecasting}\label{secF}
A possible application of suggested above filters is preliminary smoothing of
the input signals for the  predicting algorithms. For this task, the causality is crucial.
It is known that  the band-limited sequences are predictable, i.e., the sequences are predictable if with the spectrum
vanishing on a interval in $\T$.   In addition, there are predictable sequences such that the spectrum is vanishing in a single
point of $\T$; see \cite{D12sp,D12a}, where some predicable algorithms were suggested.

It can be noted that,  the suggested above filters do not change
the input sequences significantly if $a$ is close to 1; the energy of the input is not damped on a given arc  of $\T$.
In fact, the energy is damped on a small neighborhood of the point $e^{i\pi}=-1$, and the size of
this neighborhood converges to zero as $\a\to 1$. Therefore, one cannot  expect that
the filters introduced above will help to improve the performance of
the predicting algorithms  \cite{D12sp} requiring  that   the spectrum is vanishing on a fixed arc on the unit circle.

However, it appears that these filters  can help to improve  the performance  of  predictors  \cite{D12a}
oriented on processes $x$ with spectrum vanishing in a single point.
More precisely, predictors \cite{D12a} are applicable for discrete time processes $x$ such that, for some $\t>0$, $q>1$, $c>0$, \baa
\sup_{\o\in[-\pi,\pi]} |X\ew|\brea \le c
M_{\t,q}\ew,\quad X=\Z x.\label{q} \eaa
In particular, it follows that  filters (\ref{NC}) transfer sequences of a general type into predictable
sequences such that (\ref{q}) holds; respectively, filters (\ref{NC}) cannot be causal.
\par
The predicting kernel  \cite{D12a}   was defined as  $ k=
k(\cdot,\g)=\Z^{-1} K$, where \baa K(z)\defi z\left(
1-\exp\left[-\frac{\g}{z+ 1-\g^{-r} }\right]\right), \label{wK}\eaa
and where $r>0$ and $\g>0$ are  parameters.
This predictor  produces the process
\baaa
y(t)=\sum_{d=-\infty}^t k(t-d)x(t)
\eaaa
approximating $x(t+1)$ for $\g\to +\infty$ for all inputs $x$ satisfying (\ref{q}) with $q>1+2/r$.
(In the notations from \cite{D12a}, $r=2\mu/(q-1)$, where $\mu>1$, $q>1$ are the parameters ).
The function $K\ew$ approximates the function $e^{i\o}$ representing the forward one-step shift in the time domain;
the value $|K\ew-e^{i\o}|$
 is small everywhere but in a small neighborhood of $\o=\pi$. Therefore, the process  $y(t)$ represents an one step prediction of $x(t+1)$ if $X\ew$ vanishes with a certain rate at  $\o=-\pi$.  It was shown in \cite{D12a} that
\baaa
\sup_t|x(t+1)-y(t)|\to 0\quad \hbox{as}\quad \g\to 0,
\eaaa
for real sequences  $x$ such that (\ref{q}) holds, i.e., that the
prediction error vanishes as $\g\to +\infty$. Moreover, the error vanishes
  uniformly over classes of processes $x$ from some bounded sets from $\ell_{\infty}$,  such that  (\ref{q}) holds with a given $c$.

Predictors (\ref{wK}) are robust with respect to some small noise contamination, meaning that
the prediction error depends continuously on the intensity of the contaminating noise. However, for large $\g$,  the values
of $K\ew$ can be very large in a neighborhood of $\o=\pi$; in this case, the error can be
large even for a small noise.

\begin{figure}[ht]
\centerline{\epsfig{figure=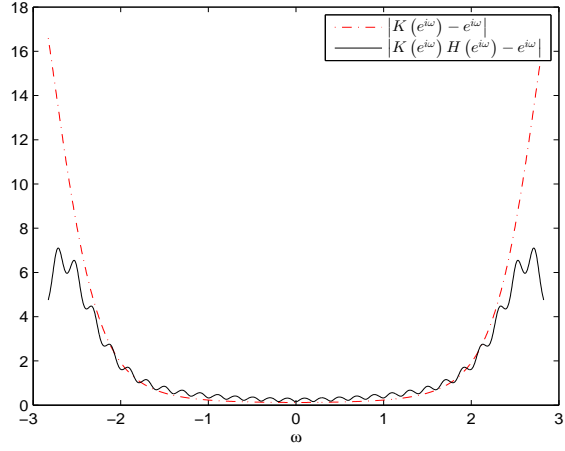,width=9cm,height=6.5cm}}
\vspace{0.5cm} \caption[]{\sm  Approximation of the one-step forward  shift operator: the values
 $|K\ew -e^{i\o}|$  for the transfer function of the predictor  (\ref{wK}) and
 $|K\ew\HH_{\nu}\ew -e^{i\o}|$,  i.e., with  smoothing
 of the input process by
filters (\ref{K}) with  $\g=2.5$, $r=0.225$, \index{[$q=10$, $\mu=1.01$,]} $\aa=0.8$, $p=0.7$, $N=30$, $m=3$. }
\vspace{0cm}
\label{figKKH}
\end{figure}

We suggest to apply filter (\ref{K}) to compensate the presence of large values of $K\ew$
in a small neighborhood of $\o=\pi$ and therefore to  reduce the impact of the presence of the high-frequency noise. This is illustrated by
 Figure \ref{figKKH} showing the shapes of error curves for approximation of
 the  forward one step shift operator. More precisely, it shows
 the shape of $|K\ew -e^{i\o}|$  for the predictor  (\ref{wK}) and
 the shape of $|K\ew\HH_{\aa}\ew -e^{i\o}|$  for the transfer functions  (\ref{K}) and (\ref{wK}), which corresponds to preliminary smoothing
 of the input process by
filters (\ref{K}). These shapes characterize imperfection of the predictors, since the transfer function $e^{i\o}$ corresponds to the one-step forward shift operator in time domain, i.e., $e^{i\o}$ represents an ideal non-causal error-free one-step ahead predictor.
 It appears that the application of the filter improves the approximation of $e^{i\o}$.
     \par

      Our setting does not involve stochastic processes and probability measure; it is oriented on smoothing the real sequences.
      However, to  provide an example of the application of  our smoothing filters,  we considered  a toy example with
      prediction of a stochastic Gaussian stationary process
       $x(t)$ evolving as an autoregression of AR(2) type
       \def\b{\beta}
       \baa
       &&x(t)=\b_1 x(t-1)+\b_2 x(t-2)+\s \eta(t), \quad t\in\ZZ.
       \label{AR2}\eaa
       Here $\eta(t)$ is a stochastic discrete time Gaussian white noise,  $
       \E \eta(t)=0$, $\E \eta(t)^2=1$. The coefficient
     $\s>0$ describes the intensity of the noise.

    For the estimation of the effectiveness of predictors, we  use the ratio
    \baa
      e(b_1,b_2)=\frac{\left(\sum_{t=1}^n|y(t-1)-x(t)|^2\right)^{1/2}}{\left(\sum_{t=1}^n|b_1 x(t-1)+b_2  x(t-2)-x(t)|^2\right)^{1/2} },\label{err}\eaa
where  $b_k\in\R$  are parameters. The values $y(t-1)$ are supposed to be the predictions, at time $t-1$, of the future values $x(t)$. Since the sequence $\{x(t)\}$ does not satisfy (\ref{q}) due to the presence of noise, a forecasting error is inevitable.
Ratio (\ref{err})  allows to compare the error of a predicting algorithm generating $y$ and the error generated by a linear predictor with the coefficients $b_1$ and $b_2$.
 More precisely,  the value $e(b_1,b_2)$ represent the ratio of the error generated by the predictor producing $y$ and the error generated with the error of the linear predictor based on the hypothesis that $\b_1=b_1$ and $\b_2=b_2$.

If the vector $(\b_1,\b_2)$ is  known, then the  optimal one step predictor of $x(t)$ is
\baa
y(t-1)=\b_1 x(t-1)+\b_2 x(t-2).
\label{ideal}\eaa
In this case, the value $n^{-1}\sum_{t=1}^n|\b_1 x(t-1)+\b_2 x(t-2)-x(t)|^2$ represents the sample mean of the squared error of this
optimal predictor with  known values of  $(\b_1,\b_2)$. Therefore, the
optimal predictor (\ref{ideal}) ensures that $e(\b_1,\b_2)\approx 1$ for a large enough $n$. Similarly, for any given $n$, an average value for $e(\b_1,\b_2)$  is also close to one for a sufficiently large number of Monte-Carlo trials, for optimal predictor (\ref{ideal}). Respectively,   any other predictor besides (\ref{ideal}), including  predictor (\ref{wK}), cannot achieve
 a lesser average value of $e(\b_1,\b_2)$ for a sufficiently large $n$  or  as an average value for a sufficiently large number of Monte-Carlo trials.

 However, in many practical situations, the value of $(\b_1,\b_2)$ is unknown, and, respectively, predictor (\ref{ideal}) cannot be used. On the other hand, predictor (\ref{wK}) does not require to know $(\b_1,\b_2)$ and can be applied  in models with unknown or random and time variable $(\b_1,\b_2)$ where predictors (\ref{ideal}) is not applicable. In other words, predictor (\ref{wK}) can be  applied  for  processes with unknown  shape
of the spectral representation.
  Therefore, it is reasonable to estimate the performance of a predictor using $e(b_1,b_2)$ with $b_k\neq\b_k$,
  for instance, with $b_k$ selected as  the expected value, or the median, or  upper or lower boundaries
  of unknown $\b_k$.

  Since it is impossible to implement convolution with infinitely supported kernels and inputs,
one has to use truncated kernels and inputs for  calculations.
In the experiments described below, we replaced  $k$ and $h_a=\Z^{-1}H_a$ by the truncated kernels
\baa
k_{d}(t)=\Ind_{\{t\le d\}}k(t),\quad h_{\aa,d}(t)=\Ind_{\{t\le d\}}h_\aa(t), \quad d>0.
\label{trun}\eaa
In other words, the original  $k=\Z^{-1}K$ and $h_a=\Z^{-1}H_a$ were used as some benchmarks; only their truncated versions
were actually implemented.

     We will use value (\ref{err}) to estimate the performance of predicting algorithms for the following two cases: \begin{itemize}
     \item The algorithm is applied without filtering and produces  $y=k_d\circ x$; we  denote by  $e_{\KK}(b_1,b_2)$ the corresponding values (\ref{err}).
     \item  The algorithm  is applied with filtering and produces  $y=(k_d\circ h_{\aa,d})\circ x$; we  denote by  $e_{\KH}(b_1,b_2)$ the corresponding values (\ref{err}).
     \end{itemize}
       In both cases, $y(t)$ is calculated
     using  historical data $\{x(s)\}_{t-d\le s\le t}$.
     \par
\begin{figure}[ht]
\label{fighV2}\centerline{\psfig{figure=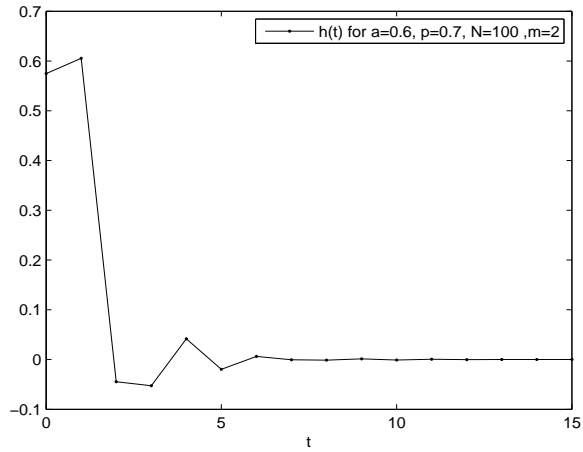, width=9cm, height=6.5cm}}
\vspace{0.5cm} \caption[]{\sm Impulse response
$h_{\a}(t)=(\F^{-1}H_{\a})(t)$ for causal filter (\ref{K})
with  the parameters given in (\ref{Par}).} \vspace{0cm}
\label{fig3}
\end{figure}
\begin{figure}[ht]
\centerline{\epsfig{figure=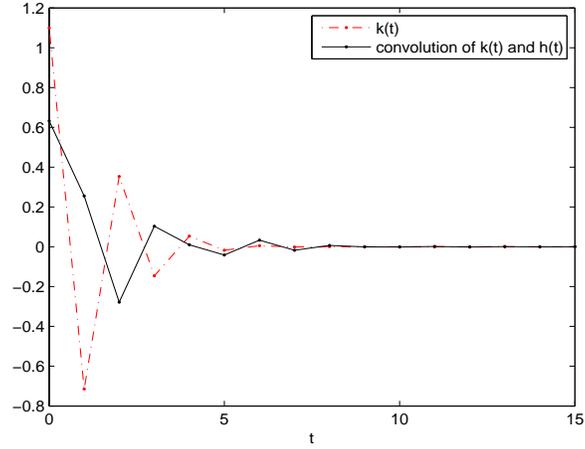, width=9cm, height=6.5cm}}
\vspace{0.5cm} \caption[]{\sm The impulse response $\Z^{-1}K$ of  predictor (\ref{wK}) and the impulse response $\Z^{-1}(K\HH_{\aa})$ of the predictor combined with
filters (\ref{K})  with the parameters given in (\ref{Par}).}
\label{figIMP}
\end{figure}

In our experiments, we used  equations (\ref{K}) and (\ref{wK}) with
 \baa
 \g=1.1,\quad r=1.1,\quad \breakk \a=0.6,\quad p=0.7,\quad N=100,\quad m=2.
 \label{Par} \eaa
 Note that selection of too large $\g$ makes  calculation of $k$ challenging, since it involves
 precise integration of fast growing $K\ew$. The choice of parameters in (\ref{Par}) ensures that the values of $|K\ew|$
 are not large.
Figure \ref{fighV2}  shows  the corresponding impulse response
$\Z^{-1}H_{\aa}$. Figure \ref{figIMP} shows  the corresponding impulse responses
$\Z^{-1}K$ and $\Z^{-1}(K\HH_{\aa})$.

     In our experiment with AR(2) process, we used  10,000
     Monte-Carlo trials with $n=d=100$ and $\s=0.3$. For each trial,  we selected   $(\b_1,\b_2)$ randomly and independently.
     The distribution of $(\b_1,\b_2)$ at each trial was the following: $\b_1$ has the uniform distribution  on the interval $(0,1)$, and $\b_2=\xi \sqrt{1-\b_1^2}$,
       where $\xi$ is a random variable  independent on $\b_1$
     and  uniformly distributed on the interval $(-1,1)$.  This choice ensures that the eigenvalues of the autoregression stay inside of the unit circle $D$ almost surely.

     We used MATLAB and standard personal computers; an experiment with 10,000 Monte-Carlo trials would take
      about five minutes of calculation time.

       First, we compared the relative
      performance of our predictors  with respect to the performance of the optimal  predictor that requires that the values $\b_1$ and $\b_2$ are known.
      We obtained that the mean value over all Monte-Carlo trials  of $e_{\KK}(\b_1,\b_2)$ is
       1.5019 and the mean value of  $e_{\KH}(\b_1,\b_2)$ is 1.1177.
     This indicates that application of filter (\ref{K}) improves the performance of the predictor.
     As was mentioned above,  we cannot expect,  for sufficiently many Monte-Carlo trials,
     the mean values of $e_{\KK}(\b_1,\b_2)$ and  $e_{\KH}(\b_1,\b_2)$ to be less than one, so the performance of the  predictor (\ref{wK})  combined with
filters (\ref{K}) is reasonably good.
\par
Second, we calculated the mean values   $e_{\KK}(b_1,b_2)$ and  $e_{\KH}(b_1,b_2)$ with $b_k=\E \b_k$,  where
$\E\b_k$ is the population mean for $\b_k$, $k=1,2$. For our parameters of the Monte-Carlo trials, we have $\E\b_1=0.5$ and $\E\b_2=0$.
 The corresponding values $e_{\KK}(0.5,0)$ and  $e_{\KH}(0.5,0)$ represent comparison of the performance of our predictor (without or with preliminary filtering)
  with an one-step predictor of $x(t)$ given by
\baaa
y(t-1)=(\E\b_1) x(t-1)+(\E\b_2) x(t-2).
\label{ideal2}\eaaa  Note that this predictor requires to know the population means of $\b_1$ and $\b_2$.
  We obtained  that the mean value for  $e_{\KK}(0.5,0)$ is $1.2292$ and the mean value for $e_{\KH}(0.5,0)$ is $0.9545$. These numbers indicate a good performance  of our
  filter/predictor system, especially if we take into account that our system does not require to know
the values $\E\b_k$.
  \par
\begin{figure}[ht]
\centerline{\epsfig{figure=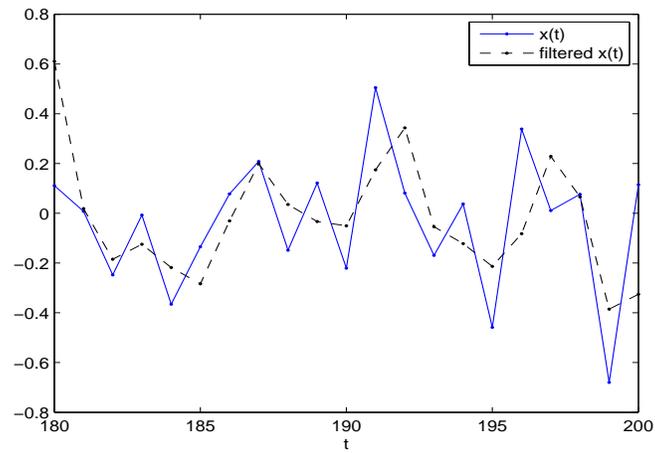,width=10cm, height=6.5cm}}
\caption[]{\sm
A path of  AR(2) process $x(t)$ versus the output of filter (\ref{K})  with the parameters given in (\ref{Par}).}
\label{figF}
\end{figure}
\begin{figure}[ht]
\centerline{\epsfig{figure=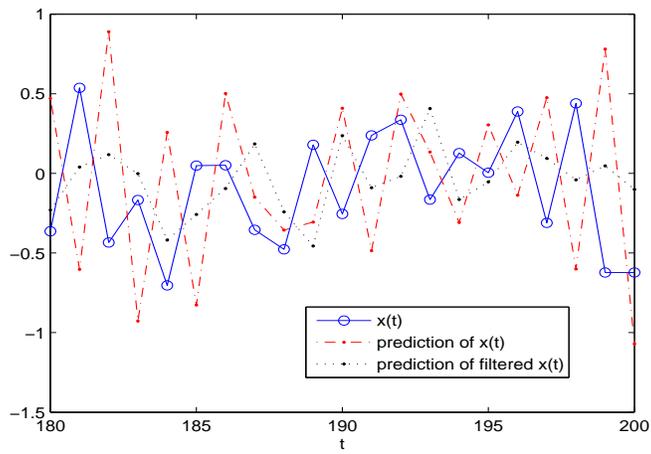,width=10cm, height=6.5cm}}
\caption[]{\sm
A path of AR(2)  process $x(t)$ versus two predictions $y(t-1)$ of $x(t)$; one was  calculated
without  filtering,  and  another was calculated after application of filter (\ref{K}) with the parameters given in (\ref{Par}).}
\label{fig4}
\end{figure}

Figure \ref{figF} shows  a sample  path of AR(2) process $x(t)$ and a filtered process  obtained using  filter (\ref{K}) with the parameters defined by (\ref{Par}). Figure \ref{fig4} shows
sample  paths of AR(2) process $x(t)$ and outputs $y(t)$
  of predictor \cite{D12a} without preliminary
    filtering and  with  preliminary filtering using filter (\ref{K})  with the parameters defined by (\ref{Par}).
     It shows the values  $y(t-1)$, i.e., predictions of $x(t)$,  versus the values of $x(t)$.

In addition, we considered a modification of process (\ref{AR2}) with $\b_2\equiv 0$, i.e., AR(1) process.
We used the same predictors and filters  as for the experiments with AR(2) process described above.

We set again  $10,000$ Monte-Carlo trials with $n=d=100$, with $\s=0.3$, and with randomly selected  $\b_1$
     such that $\b_1$ was distributed uniformly on the interval $(0,1)$.
     We compared the relative
      performance of our predictors  with respect to the performance of the optimal  predictor that requires that the value $\b_1$ is known.
      We obtained that the mean value of $e_{\KK}(\b_1,0)$ is
       1.2830 and the mean value of  $e_{\KH}(\b_1,0)$ is 1.1023.
     The numbers  indicate  again that the use of the filter improves the performance of the predictor.
     Again, the mean values of $e_{\KK}(\b_1,0)$ and  $e_{\KH}(\b_1,0)$ cannot be less than one, so the performance of the  predictor (\ref{wK})  combined with
filters (\ref{K}) is reasonably good.

Finally, we calculated the values   $e_{\KK}(b_1,0)$ and  $e_{\KH}(b_1,0)$ with $b_1=\E \b_1=0.5$ for AR(1) process defined by (\ref{AR2}) with $\b_2\equiv 0$. In other words, we compared the relative
performance of the one step predictor (\ref{K})  with or without filtering
with respect to the predictor $y(t)=(\E\b_1) x(t)$  that requires to known the population mean for $\b_1$. We obtained  that  $e_{\KK}(0.5,0)= 1.1824$ and  $e_{\KH}(0.5,0)=1.0155$. This shows again that application of the filter reduces the forecasting error.  Given that our system does not require to know  $\E\b_1$, the performance is reasonably good.

 The results of these experiments appear to be consistent and stable  with respect to the variations of the parameters.
\subsection*{On the impact of truncation}
Since it is impossible to implement convolution with infinitely supported kernels and inputs,
we have to  run numerical calculations with truncated processes.  Let us show  that the described above
filtering and  forecasting  are robust with respect to the truncation, i.e.,  that  the truncation
 in  (\ref{trun}) has a vanishing impact  for large $d$.  Since $H_\aa\ew\in L_\infty(-\pi,\pi)$, we have that $h_\aa\in\ell_2$. Hence $h_{\aa,d}\circ x\to h_\aa\circ x$  in $\ell_\infty$ as $d\to +\infty$ for $x\in\ell_2$; in practice, only truncated inputs $x\in \ell_2$ are available. It can be also
noted that the predicting kernel $k=\Z^{-1}K$ defined by (7)  belongs to $\ell_2$. Therefore, the kernel $k\circ h_{\aa,d}$
converges to $k\circ h_{\aa}$ in $\ell_\infty$ as $d\to +\infty$.

Our  numerical experiments for autoregressions with the coefficients deemed to untraceable
demonstrated that truncation
with relatively small $d=100$ does not diminish a good forecasting performance.
 \section{Conclusion}
The paper proposes a family of causal smoothing filters. These
filters are near-ideal meaning that a higher rate of damping of the energy on the
high frequencies    would lead to the loss of causality; this is because they
approximate non-causal filters transferring non-predicable processes into predictable ones.  A
possible application is preliminary smoothing of the inputs for predicting algorithms. Certain mild but
stable improvement of forecasting accuracy is demonstrated in experiments with simple autoregressions in a setting where theirs coefficients are deemed to be untraceable.

  It could be interesting to investigate the computational limits of the algorithms described above,
for instance, for long sequences, for higher order autoregressions, or for other types of input processes.
It could be interesting to find other filters with similar properties. It could be useful
to represent the filtering algorithm in the terms of the discrete Fourier transform. We leave this  for future work.

   \comm{In a more general setting, it could be also interesting to consider minimization of the error
 $\|\w K\ew -e^{i\o}\|_{L_2(-\pi,\pi)}$ among all functions  $\w K\ew$ from a Hardy space on $D^c$.}
\subsection*{Acknowledgment}
This work  was supported by ARC grant of Australia DP120100928 to the author.
In addition, the author thanks the anonymous referees for their valuable remarks which  helped to improve the paper.



\end{document}